\documentclass[11pt]{article}
\usepackage{amsfonts,amssymb}
\usepackage{amsmath,amscd}
\begin{document}
\sloppy
\newcommand{\dickebox}{{\vrule height5pt width5pt depth0pt}}
\newtheorem{Def}{Definition}[section]
\newtheorem{Bsp}{Example}[section]
\newtheorem{Prop}[Def]{Proposition}
\newtheorem{Theo}[Def]{Theorem}
\newtheorem{Lem}[Def]{Lemma}
\newtheorem{Koro}[Def]{Corollary}
\newcommand{\Hom}{{\rm Hom}}
\newcommand{\End}{{\rm End}}
\newcommand{\Ext}{{\rm Ext}}
\newcommand{\add}{{\rm add}}
\newcommand{\rad}{{\rm rad}}
 \newcommand{\lra}{\longrightarrow}
 \newcommand{\ra}{\rightarrow}
{\Large \bf
\begin{center} Transfer of stable equivalences of Morita type
 \end{center}}
\medskip

\centerline{{\sc Shengyong Pan} and {\sc Changchang Xi $^*$}}
\begin{center} School of Mathematical Sciences, \\ Laboratory of Mathematics and Complex Systems,\\
Beijing Normal University, 100875 Beijing,\\
People's Republic of  China \\
\texttt{E-mail:xicc@bnu.edu.cn}\\
\texttt{E-mail:panshy1979@bnu.edu.cn}\\
\end{center}
\renewcommand{\thefootnote}{\alph{footnote}}
\setcounter{footnote}{-1} \footnote{ $^*$ Corresponding author.
Email: xicc@bnu.edu.cn; Fax: 0086 10 58802136; Tel.: 0086 10
58808877.}

\medskip
\abstract{Let $A$ and $B$ be finite-dimensional $k$-algebras over a
field $k$ such that $A/\rad(A)$ and $B/\rad(B)$ are separable. In
this note, we consider how to transfer a stable equivalence of
Morita type between $A$ and $B$ to that between $eAe$ and $fBf$,
where $e$ and $f$ are idempotent elements in $A$ and in $B$,
respectively. In particular, if the Auslander algebras of two
representation-finite algebras $A$ and $B$ are stably equivalent of
Morita type, then  $A$ and $B$ themselves are stably equivalent of
Morita type.  Thus, combining a result with Liu and Xi, we see that
two representation-finite algebras $A$ and $B$ over a perfect field
are stably equivalent of Morita type if and only if their Auslander
algebras are stably equivalent of Morita type. Moreover, since
stable equivalence of Morita type preserves $n$-cluster tilting
modules, we extend this result to $n$-representation-finite algebras
and $n$-Auslander algebras studied by Iyama.}

\medskip {\small {\it 2000 AMS Classification}: 16G10, 18G05,
16S50; 16P10, 20C05, 18E10.

\medskip {\it Key words:} stable equivalence of Morita type,
Auslander algebra, $n$-cluster tilting module.}

\section{Introduction}
Stable equivalence of Morita type introduced by Brou\'{e} in
\cite{B} is one of the fundamental equivalence relations for
algebras and groups. It is of considerable interest in the modular
representation theory of finite groups, or more generally, of
finite-dimensional self-injective algebras. The notion is intimately
related to the celebrated conjecture of Brou\'{e}, which says that
certain blocks of group algebras with abelian defect groups should
be derived-equivalent (see \cite{B} and \cite{R1} for precise
formulation); the connection can be seen from Rickard's result that
a derived equivalence between self-injective algebras induces a
stable equivalence of Morita type \cite{R2}. Recently, it is shown
that stable equivalence of Morita type is also of particular
interest for general finite-dimensional algebras, for example, it
preserves representation dimension \cite{C2}, representation type
\cite{krause}, Hochschild homological and cohomological groups
\cite{LX2, C3}, and the absolute value of Cartan determinant
\cite{C3}. As is known, stable equivalence of Morita type occurs
frequently not only in the block theory of finite groups \cite{L2},
but also in the representation theory of general finite-dimensional
algebras. A plenty of such examples are constructed in
\cite{LX1,LX2,LX3}.

Moreover, it is shown in \cite{LX2} that, for two finite-dimensional
$k$-algebras $A$ and $B$ over a field $k$ of finite
representation-type, if $A$ and $B$ are stably equivalent of Morita
type, then their Auslander algebras are also stably equivalent of
Morita type. A natural question is whether the converse of this
statement is true.

In this note, we shall consider the general question of how to
transfer a stable equivalence of Morita type between two algebras
$A$ and $B$ over a field to that between $eAe$ and $fBf$, where $e$
and $f$ are idempotent elements in $A$ and $B$, respectively. We say
that two bimodules $_AM_B$ and $_BN_A$ define a stable equivalence
of Morita type between $A$ and $B$ if $M$ and $N$ are projective as
one-sided modules, and there are a projective $A$-$A$-bimodule $P$
and a projective $B$-$B$-bimodule $Q$ such that $M\otimes_BN\simeq
A\oplus P$ and $N\otimes_AM\simeq B\oplus Q$ as bimodules. With
these notations in mind, our main result reads as follows:

\begin{Theo} Suppose that $A$ and $B$ are finite-dimensional $k$-algebras over
a field $k$ such that both $A/\rad(A)$ and $B/\rad(B)$ are
separable. Let $_AM_B$ and $_BN_A$ be two bimodules defining a
stable equivalence of Morita type between $A$ and $B$. If $e^2=e\in
A$ such that $Pe\in \add(Ae)$, and if $f$ is a sum of the pairwise
orthogonal primitive idempotent elements $f_j\in B$ in the
decomposition $Ne\simeq \oplus_{j=1}^m (Bf_j)^{r_j}$, then the
bimodules $eMf$ and $fNe$ define a stable equivalence of Morita type
between $eAe$ and $fBf$. \label{thm}
\end{Theo}

From this result, we have the following corollary which supplies a
positive answer to our previous question on Auslander algebras. For
the unexplained notion in the corollary, we refer the reader to
Section \ref{sec3} below.

\begin {Koro} Suppose that $A$ and $B$ are finite-dimensional $k$-algebras over a perfect field $k$.
Assume that $A$ and $B$ are $n$-representation-finite. If the
$n$-Auslander algebras of $A$ and $B$ are stably equivalent of
Morita type, then $A$ and $B$ themselves are stably equivalent of
Morita type.
\end{Koro}

The proof of Theorem \ref{thm} is given in the next section.

\section{Proof of the main result \label{sec2}}

Throughout this note, $k$ denotes a fixed field. Given a
finite-dimensional $k$-algebra $A$, we denote by $A$-mod the
category of all finitely generated left $A$-modules. If $M\in
A$-mod, we denote by add$(M)$ the full subcategory of $A$-mod
consisting of all modules $X$ which are direct summands of finite
sums of copies of $M$. By an algebra we mean a finite-dimensional
$k$-algebra, and by a module we mean a left module. The global and
dominant dimensions of an algebra $A$ are denoted by gl.dim$(A)$ and
dom.dim$(A)$, respectively. The composition of two homomorphisms
$f:X\ra Y$ and $g: Y\ra Z$ is denoted by $fg: X\ra Z$, and the usual
$k$-duality is denoted by $D:=\Hom_k(-, k)$. Let us recall the
definition of a stable equivalence of Morita type.

\begin{Def}$\cite{B}$ Suppose that $A$ and $B$ are two $($arbitrary$)$ $k$-algebras.
We say that $A$ and $B$ are stably equivalent of Morita type if
there exist an $A$-$B$-bimodule $_{A}M_{B}$ and a $B$-$A$-bimodule
$_{B}N_{A}$ such that

$(1)$ $M$ and $N$ are projective as one-sided modules, and

$(2)$ $M\otimes N\simeq A\oplus P $ as $A$-$A$-bimodules for some
projective $A$-$A$-bimodule $P$,and $N\otimes M\simeq B\oplus Q $ as
$B$-$B$-bimodules for some projective $B$-$B$-bimodule $Q$.
\end {Def}

Note that if $A$ and $B$ are stably equivalent of Morita type, then
their opposite algebras $A^{op}$ and $B^{op}$ are also stably
equivalent of Morita type.

Let $A$ be a representation-finite algebra. An $A$-module $X$ is
called an additive generator for $A$-mod if add($X$) = $A$-mod, that
is, every indecomposable $A$-module is isomorphic to a direct
summand of $X$. Let $X$ be an additive generator for $A$-mod. The
endomorphism algebra $\Lambda= \End_{A}(X)$ of $X$ is called the
Auslander algebra of $A$. (This is unique up to Morita equivalence.)
Auslander algebras can be described by two homological properties:
An algebra $A$ is an Auslander algebra if (1) gl.dim($A$)$\leq 2$;
and (2) if $0\ra A \ra I_{0}\ra I_1\ra I_2\ra 0$ is a minimal
injective resolution of $A$, then $I_{0}$ and $I_{1}$ are
projective.

An $A$-module $X\in A$-mod is called a generator for $A$-mod if
$\add(_AA)\subseteq \add(X)$; a cogenerator for $A$-mod if
$\add(D(A_A))\subseteq \add(X)$, and a generator-cogenerator if it
is both a generator and a cogenerator for $A$-mod. Clearly, for a
representation-finite algebra $A$, an additive generator for $A$-mod
is a generator-cogenerator for $A$-mod.

In the following, we shall introduce some notations. Assume that $A$
is a $k$-algebra.

Let $T$ be an arbitrary $A$-module, and suppose $B$ is the
endomorphism algebra of $T$. We consider the following subcategories
related to $T$.

$Gen(_{A}T)$ =$\{X\in A$-mod $\mid$ there is a surjective
homomorphism from $T^{m}$ to $X$ with $m\ge 1 \}$.

$Pre(_{A}T)$ =$\{X\in A$-mod $\mid$ there is an exact sequence
$T_{1}\rightarrow T_{0}\rightarrow X $ with all $T_i\in$
add$(_{A}T)$ \}.

$App(_AT)$ =$\{X\in A$-mod $\mid $ there is a homomorphism $g:
T_0\ra X$ with  $T_0\in$ add$(_AT)$ such that
\\ Ker$(g)\in Gen(_{A}T)$ and $g$ is a right add$(_{A}T)$-approximation
of $X\}$.

\medskip
The following lemma is known, for a proof, we refer to \cite[Lemma
2.1]{C5}.

\begin{Lem} Let $X$ be an arbitrary $A$-module. Recall that
$B=\End_A(T)$ and $_AT_B$ is a natural bimodule. Then:

$(1)$ Let $Y$ be a right $B$-module. The natural homomorphism
$\delta$: $Y\otimes_B \Hom_{A}(T,X)\rightarrow
 \Hom_{B}(\Hom_{A}(X,T), Y)$, given by $y\otimes f \mapsto \delta_{y\otimes f}$
with $\delta_{y\otimes f}(g)= y(fg)$ for $y\in Y, f\in \Hom_A(T,X),
g\in \Hom_A(X,T)$, is an isomorphism if $X\in
 \emph{add}(_{A}T)$.

$(2)$ If $X'\in \add (_AT)$, or $ X \in \add(_AT)$, then the
composition map $\mu: \Hom_A(X',T)\otimes_B\Hom_A(T,X)\rightarrow
\Hom_A(X',X)$ given by $f\otimes_B g\mapsto fg $ is bijective.

$(3)$ Let $C$ be a $k$-algebra, and suppose $_AX_C$ is an
$A$-$C$-bimodule. If $_AX\in Gen(_{A}T)$, then the evaluation map
$e_{X}: T\otimes_{B}\Hom_A(T,X)\rightarrow X$ is surjective as
$A$-$C$-bimodules. If $X\in App(_{A}T)$, then $e_{X}$ is an
isomorphism as $A$-$C$-bimodules. Conversely, if $e_{X}$ is
bijective as $A$-modules, then $X\in App(_{A}T)$. \label{2.3}
\end{Lem}

The next lemma is taken from \cite[Lemma 2.1]{C3}.

\begin{Lem}$ \cite{C3}$ $(1)$ Let $A, B, C$
and $E$ be $k$-algebras, and let $_{A}X_{B}$ and $_{B}Y_{E}$ be
bimodules with $X_{B}$ projective. Put $X^{*}=\Hom_{B}(X,B)$. Then
the natural homomorphism $\phi: {}_AX\otimes_{B}Y_{E}\rightarrow
\Hom_{B}(_BX^*_A, {}_BY_E)$, defined by $f\mapsto (xf)y$ for $x\in
X, y\in Y$ and $f\in X^{*}$, is an isomorphism of $A$-$E$-bimodules,
where the image of $x$ under $f$ is denoted by $xf$.

$(2)$ In the situation $(_EP_A, {}_CX_{B}, {}_AU_B)$, if $P_A$ is
projective, or if $X_B$ is projective, then
$_EP\otimes_A\Hom_B(_CX_B, {}_AU_B)\simeq \Hom_B(_CX_B,
{}_EP\otimes_{A}U_B)$ as $E$-$C$-bimodules. Dually, in the situation
$(_AP_E, {}_BX_C, {}_BU_A)$, if $ _AP$ is projective, or if $_{B}X$
is projective, then $\Hom_B(_BX_C, {}_BU_A)\otimes_AP_E\simeq
\Hom_B(_BX_C, {}_BU\otimes_AP_E)$ as $C$-$E$-bimodules.\label{2.4}
\end{Lem}

The following is a well-known result due to Auslander (for example,
see \cite[Proposition 5.6, p.214]{ARS}).

\begin{Lem} Let $\Lambda$ be an Artin algebra such that $\emph{gl.dim}(\Lambda)\leq 2\leq
 \emph{dom.dim}(\Lambda)$. Let $U$ be a $\Lambda$-module such that
$\add(U)$ is the full subcategory of $\Lambda\emph{-mod}$ consisting
of all projective-injective $\Lambda$-modules. Then

$(1)$ $A:= \End_{\Lambda}(U)$ is representation-finite.

$(2)$ $\Lambda$ is Morita equivalent to $\End_A(X)$, where $X$ is an
additive generator for $A\emph{-mod}$. \label{2.5}
\end{Lem}

For our proof of Theorem \ref{thm}, we also need the following lemma
in \cite[Theorem 2.7, Corollary 3.1, Lemma 3.2]{DM}, see also
\cite[Lemma 3.3]{C3}.

\begin{Lem}$\cite{DM}$ Suppose that $A$ and $B$ are finite-dimensional $k$-algebras over a field $k$
such that $A$ and $B$ have no separable direct summands and that
$A/\rad(A)$ and $B/\rad(B)$ are separable. Assume that $_AM_B$ and
$_BN_A$ are indecomposable bimodules that define a stable
equivalence of Morita type between $A$ and $B$. Then:

$(1)$ There are isomorphisms of bimodule: $N\cong \Hom_{A}(M,A)\cong
\Hom_{B}(M,B)$ and $M\cong \Hom_{A}(N,A)\cong \Hom_{B}(N,B)$.

$(2)$ Both $(N \otimes_{A} -, M \otimes_{B} - )$ and $( M
\otimes_{B} -, N\otimes_{A}-)$ are adjoint pairs of functors.

$(3)$ If $_{A}I$ is injective, then so is $N\otimes_A I$.
 \label{2.6}
\end{Lem}

It follows from Lemma \ref{2.6} that the following result is true.
Note that the last statement in Lemma \ref{2.8} below follows from
\cite[Lemma 4.5]{C3}.

\begin{Lem} Suppose that $A$ and $B$ are finite-dimensional $k$-algebras over a field $k$
such that $A/\rad(A)$ and $B/\rad(B)$ are separable. Assume that
$\{e_1,\cdots, e_n\}$ and $\{f_1, \cdots, f_m\}$ are complete sets
of pairwise orthogonal primitive idempotents in $A$ and $B$,
respectively. Let $e$ be the sum of all those $e_i$ for which $Ae_i$
is projective-injective, and let $f$ be the sum of all those $f_j$
for which $Bf_j$ is projective-injective. If $M$ and $N$ are
indecomposable bimodules that define a stable equivalence of Morita
type between $A$ and $B$, then $Ne\simeq N\otimes_{A}Ae \in
\add(Bf), Mf\simeq M\otimes_{B}Bf \in \add(Ae)$, and $Pe\in
\add(Ae)$. \label{2.8}
\end{Lem}

{\bf Proof of Theorem 1.1:}

Suppose that $A$ and $B$ are two algebras over a  field $k$ such
that $A/\rad(A)$ and $B/\rad(B)$ are separable, and suppose that
$_AM_B$ and $_BN_A$ define a stable equivalence of Morita type
between $A$ and $B$. We may assume that both $A$ and $B$ have no
separable summands since the direct sum of $A$ with a separable
$k$-algebra is always stably equivalent of Morita type to $A$
itself. Furthermore, by \cite[Lemma 2.2]{LX1}, we may assume that
$M$ and $N$ are indecomposable as bimodules. Then $(M\otimes_B-,
N\otimes_A-)$ and $(N\otimes_A-, M\otimes_B-)$ are adjoint pairs by
Lemma \ref{2.6}.

To prove Theorem \ref{thm}, we shall show that the bimodules $eMf$
and $fNe$ satisfy the conditions of a stable equivalence of Morita
type between $eAe$ and $fBf$.

$(1)$ $fNe$ is projective as both an $fBf$-module and a right
$eAe$-module.

In fact, we have $fNe \simeq fB\otimes_BNe \simeq \Hom_B(Bf,
B)\otimes_BNe \simeq \Hom_B(Bf, {}_BNe)$ by Lemma \ref{2.4}. Since
$Ne \in \add(Bf)$ by the definition of $f$, we see that
$\Hom_{B}(Bf, Ne)$ is projective as an $fBf$-module, that is, $fNe$
is projective as an $fBf$-module. To see that $fNe$ is a projective
right $eAe$-module, we notice that
$\add(Mf)=\add(M\otimes_BBf)=\add(M\otimes_BNe)=\add((M\otimes_BN)e)$
= $\add(Ae\oplus Pe)=\add(Ae)$, here we use the assumption $Pe\in
\add(Ae)$. Since $(M \otimes_B -, N \otimes_A - )$ is an adjoint
pair, it follows from $\Hom_{B}(Bf, {}_{B}N\otimes_{A}Ae)\simeq
\Hom_{A}(M\otimes_B Bf, Ae)\simeq \Hom_A(Mf,Ae)$ that $fNe$ is
projective as a right $eAe$-module since $Mf \in \add(Ae)$. Thus (1)
is proved.

$(2)$ $eMf$ is projective as both an $eAe$-module and a right
$fAf$-module. The proof of (2) is similar to that of (1), we omit it
here.

$(3)$ $fNe\otimes eMf \simeq fBf \oplus fQf$ as bimodules.

Indeed, by the associativity of tensor products, we have the
following isomorphisms of $fBf$-$fBf$-bimodules:

$$\begin{array}{rl} fNe\otimes_{eAe} eMf & \simeq fN\otimes_{A}Ae\otimes_{eAe}eA\otimes_{A}Mf \\
& \simeq fN\otimes_{A}Ae\otimes_{eAe}\Hom(Ae, A)\otimes_{A}Mf\\
& \simeq fN\otimes_AAe\otimes_{eAe}\Hom(Ae, {}_AMf) \quad (\mbox{ by
Lemma \ref{2.4}}\;)\\
& \simeq fN\otimes_AMf \quad (\mbox{by Lemma \ref{2.3}}\;).
\end{array}$$
Since $M$ and $N$ define the stable equivalence of Morita type
between $A$ and $B$, we have $N\otimes_A M\simeq B\oplus Q$ as
$B$-$B$-bimodules. This implies that $fNe\otimes eMf \simeq
fN\otimes Mf \simeq fB \otimes_{B}N \otimes_A M \otimes_{B} Bf\simeq
fB\otimes_{B} (B\oplus Q) \otimes_{B} Bf \simeq fBf \oplus fB
\otimes_B Q \otimes_B Bf \simeq fBf \oplus fQf$.

$(4)$ The bimodule $fQf$ in (3) is projective.

In fact, since $Q$ is a projective $B$-$B$-bimodule and since
$B/\rad(B)$ is separable, the bimodule $Q$ is isomorphic to a direct
sum of modules of the form $Q_{1}\otimes_k Q_{2}$, where $Q_{1}$ is
a projective left $B$-module and $Q_{2}$ is a  projective right
$B$-module. Since $fNe\otimes_{eAe} eMf$ and $fBf$ are projective as
left $fBf$-modules, we infer that $fQ_{1}\otimes_{k} Q_{2}f$ is a
projective $fBf$-module. It follows that $fQ_{1}$ is a projective
$fBf$-module. Similarly, $Q_{2}f$ is a  projective right
$fBf$-module. Thus $fQf$, which is isomorphic to a direct sum of
modules of the form $fQ_{1}\otimes_{k} Q_{2}f$ with $fQ_{1}$ a
projective $fBf$-module and $Q_{2}f$ a projective right
$fBf$-module, is projective as an $fBf$-$fBf$-bimodule.

$(5)$ Similarly, we can show that $eMf\otimes_{fBf}fNe\simeq
eAe\oplus ePe$ and that $ePe$ is a projective $eAe$-$eAe$-bimodule.

Indeed, we have
$$\begin{array}{rl} eMf\otimes_{fBf} fNe & \simeq eM\otimes_BBf\otimes_{fBf}fB\otimes_{B}Ne \\
& \simeq eM\otimes_{B}Bf\otimes_{fBf}\Hom(Bf, B)\otimes_{B}Ne\\
& \simeq eM\otimes_BBf\otimes_{fBf}\Hom(Bf, {}_BNe) \quad (\mbox{ by
Lemma \ref{2.4}}\;)\\
& \simeq eM\otimes_BNe \quad (\mbox{by Lemma \ref{2.3}}\;).
\end{array}$$
Since $M$ and $N$ define the stable equivalence of Morita type
between $A$ and $B$, we have $M\otimes_BN\simeq A\oplus P$ as
$A$-$A$-bimodules. This implies that $eMf\otimes_{fBf} fNe \simeq
eM\otimes_ANe \simeq eAe \oplus ePe$. Now, we show that the bimodule
$ePe$ is projective.

In fact, since $P$ is a projective $A$-$A$-bimodule and since
$A/\rad(A)$ is separable, the bimodule $P$ is isomorphic to a direct
sum of modules of the form $P_{1}\otimes_k P_{2}$, where $P_{1}$ is
a  projective left $A$-module and $P_{2}$ is a  projective right
$A$-module. Since $eMf\otimes_{fBf} fNe$ and $eAe$ are projective as
left $eAe$-modules, we infer that $eP_{1}\otimes_{k} P_{2}e$ is a
projective $eAe$-module. It follows that $eP_1$ is a projective
$eAe$-module. Similarly, $P_{2}e$ is a projective right
$eAe$-module. Thus $ePe$, which is isomorphic to a direct sum of
modules of the form $eP_{1}\otimes_{k} P_{2}e$ with $eP_{1}$ a
projective $eAe$-module and $P_{2}e$ a projective right
$eAe$-module, is projective as an $eAe$-$eAe$-bimodule.

Thus, by definition, the bimodules $eMf$ and $fNe$ define a stable
equivalence of Morita type between $eAe$ and $fBf$. This finishes
the proof of Theorem \ref{thm}. $\square$

\medskip
{\it Remarks}. (1) In Theorem \ref{thm}, if $e$ is an idempotent
element in $A$ such that every indecomposable projective-injective
$A$-module is isomorphic to a summand of $Ae$, then $Pe\in
\add(Ae)$. This follows immediately from the proof of \cite[Lemma
4.5]{C3}. In fact, under the assumption of Theorem \ref{thm}, we
infer that $\nu^i_A(_AP)$ is projective-injective for all $i\ge 0$,
where $\nu_A$ is the Nakayama functor $D\Hom_A(-, {}_AA)$. Hence, if
$e$ is an idempotent element in $A$ such that every indecomposable
projective-injective $A$-module $X$ with $\nu^j_AX$
projective-injective for all $j\ge 0$ belongs to $\add(Ae)$, then
$_AP\in \add(Ae)$.

(2) As was pointed out in \cite[Section 4]{DM}, if $e$ is an
idempotent in $A$ and if $f$ is an idempotent in $B$ such that
$\add(Ae)$ and $\add(Bf)$ are invariant under Nakayama functor, then
$eAe$ and $fBf$ are self-injective, and any stable equivalence of
Morita type between $A$ and $B$ induces a stable equivalence of
Morita type between $eAe$ and $fBf$. In general, however,  our
algebras $eAe$ and $fBf$ in Theorem \ref{thm} may not be
self-injective.

\medskip
As a corollary of Theorem \ref{thm}, we get the following result.

\begin {Koro}
Suppose $A$ and $B$ be are finite-dimensional $k$-algebras of finite
representation type, and let $\Lambda$ and $\Gamma$ be the
corresponding Auslander algebras of $A$ and $B$, respectively.
Assume that $\Lambda/\rad(\Lambda)$ and $\Gamma/\rad(\Gamma)$ are
separable. Then $\Lambda$ and $\Gamma$ are stably equivalent of
Morita type if and only if A and B are stably equivalent of Morita
type. \label{cor3}
\end{Koro}

{\it Proof.} We know that if $X$ is an additive generator for
$A$-mod with $\Lambda := \End_A(X)$, then $U:=\Hom_A(X, D(A_A))$ is
a projective-injective $\Lambda$-module with
$\End_{\Lambda}(U)\simeq A^{op}$; and every indecomposable
projective-injective $\Lambda$-module is isomorphic to a direct
summand of $U$. Note that this $U$ satisfies the conditions in Lemma
\ref{2.5}. If we choose $e$ to be the sum of all idempotents
corresponding to the indecomposable injective $A$-modules, then
Lemma \ref{2.8} says that the conditions in Theorem \ref{thm} on the
idempotent $e\in \Lambda$ are satisfied. Note that $e$ defines an
idempotent element $f$ in $\Gamma$ (see Theorem \ref{thm}), and that
$\add(\Gamma f)$ contains all projective-injective $\Gamma$-modules.
With these in mind, the corollary follows from Theorem \ref{thm} and
\cite[Theorem 1.1]{LX3}. $\square$

\medskip
For an algebra $A$, we denote by $[A]$ the class of all those
algebras $B$ for which there is a stable equivalence of Morita type
between $B$ and $A$. From the above corollary, we have the following
result.

\begin {Koro} Suppose that $k$ is a perfect field. Let $\cal F$ be the set of equivalence
classes $[A]$ of representation-finite $k$-algebras $A$ with respect
to stable equivalence of Morita type, and let $\cal A$ be the set of
equivalence classes $[\Lambda]$ of Auslander $k$-algebras $\Lambda$
with respect to stable equivalence of Morita type. Then there is an
one-to-one correspondence between $\cal F$ and $\cal A$.
\end{Koro}

Another consequence of Theorem \ref{thm} is the following corollary.

\begin{Koro} Suppose that $A$ and $B$ are two $k$-algebras.
Let $_AX$ be a generator-cogenerator for $A\emph{-mod}$ such that
$\End_A(X)/\rad(\End_A(X))$ is separable, and let $_BY$ be a
generator-cogenerator for $B\emph{-mod}$ such that
$\End_B(Y)/\rad(\End_B(Y))$ is separable. If $\End_A(X)$ and
$\End_B(Y)$ are stably equivalent of Morita type, then so are $A$
and $B$.  In this case, $A$ and $B$ have the same global, dominant,
finitistic and representation dimensions.
\end{Koro}

Finally, we remark that if we consider derived equivalence instead
of stable equivalence of Morita type in Corollary \ref{cor3}, then
we know from \cite{hx4} that a derived equivalence between
representation-finite, self-injective algebras $A$ and $B$ implies a
derived equivalence between their Auslander algebras. But the
converse of this statement is still open. For further information on
constructing derived equivalences, we refer the reader to a current
paper \cite{hx2}.

\section{Higher Auslander algebras \label{sec3}}
In the following, we point out that Corollary \ref{cor3} holds true
for $n$-representation-finite algebras and $n$-Auslander algebras
studied in \cite{iyama}.

Now we recall some definitions from \cite{iyama}. Let $A$ be a
finite-dimensional $k$-algebra, and let $n\ge 1$ be a natural
number. An $A$-module $T$ is called an $n$-cluster tilting module if
$\add(T)=\{X\in A$-mod $\mid \Ext^i_A(X,T)=0, 1\le i< n\}$ = $\{X\in
A$-mod $\mid \Ext^i_A(T,X)=0, 1\le i<n\}$. The $k$-algebra $A$ is
called $1$-representation-finite if there is an $1$-cluster tilting
$A$-module $T$. This is equivalent to saying that $A$ is
representation-finite. For $n\ge 2$, the $k$-algebra $A$ is called
$n$-representation-finite if gl.dim$(A)\le n$ and there is an
$n$-cluster tilting $A$-module $T$.

A $k$-algebra $\Lambda$ is called an $n$-Auslander algebra if there
is an $n$-representation-finite $k$-algebra $A$ with an $n$-cluster
tilting $A$-module $T$ such that $\Lambda$ is Morita equivalent to
$\End_A(T)$. Note that, for an $n$-representation-finite algebra
$A$, its $n$-Auslander algebra is unique up to Morita equivalence.

Clearly, each $n$-cluster tilting $A$-module $T$ is a generator and
co-generator for $A$-mod. Thus the indecomposable
projective-injective $\End_A(T)$-modules are of the form
$\Hom_A(T,I)$, where $I$ is an indecomposable injective $A$-module.

Let $A$ be an $n$-representation-finite $k$-algebra with $T$ an
$n$-cluster tilting $A$-module. Furthermore, we assume that $A$ has
no separable direct summands and that $A/\rad(A)$ is separable. If
$A$ is stably equivalent of Morita type to an algebra $B$ such that
$B$ has no separable direct summand and $B/\rad(B)$ is separable,
then $B$ is $n$-representation-finite. In fact, if two
indecomposable bimodules $_AM_B$ and $_BN_A$ define the stable
equivalence of Morita type between $A$ and $B$, then $N\otimes_AT$
is an $n$-cluster tilting $B$-module: Since this stable equivalence
of Morita type is of adjoint type by Lemma \ref{2.6}, we see that
$\Ext^i_A(N\otimes_BT, N\otimes_AT)\simeq \Ext^i_A(T, M\otimes_B
N\otimes_AT)\simeq \Ext^i_A(T,T\oplus P\otimes_AT)=
\Ext^i_A(T,T)\oplus \Ext^i_A(T, P\otimes_AT)=0$ for $1\le i <n$
since $P\otimes_AT$ is a projective-injective $A$-module. This shows
that $\add(N\otimes_AT)$ is contained  in both $\{X\in B$-mod $\mid
\Ext^i_B(X,N\otimes_AT)=0, 1\le i< n\}$ and $\{X\in B$-mod $\mid
\Ext^i_B(N\otimes_AT,X)=0, 1\le i< n\}$. Now, let $Y\in B$-mod such
that $\Ext^j_B(N\otimes_AT,Y)=0$ for $1\le j<n$. Then
$0=\Ext^j_B(N\otimes_AT,Y)=\Ext^j_A(T,M\otimes_BY)$ for $1\le j<n$,
and therefore $M\otimes_BY\in \add(T)$. This implies that $Y\in
\add(N\otimes_AT)$. Similarly, we show that
$\add(N\otimes_AT)=\{Y\in B$-mod $\mid \Ext^i_B(Y,N\otimes_AT)=0,
1\le i<n \}$. Hence $N\otimes_AT$ is an $n$-cluster tilting
$B$-module.

Thus, $n$-representation-finite $k$-algebras $A$ with $A/\rad(A)$
separable are closed under stable equivalences of Morita type.

As in the case of Corollary \ref{cor3}, the following is a
consequence of Theorem \ref{thm}.

\begin {Theo}
Suppose that $A$ and $B$ are finite-dimensional $k$-algebras such
that both are $n$-representation-finite. Let $\Lambda$ and $\Gamma$
be the corresponding $n$-Auslander algebras of $A$ and $B$,
respectively. Assume that both $\Lambda/\rad(\Lambda)$ and
$\Gamma/\rad(\Gamma)$ are separable. Then $\Lambda$ and $\Gamma$ are
stably equivalent of Morita type if and only if $A$ and $B$ are
stably equivalent of Morita type.\label{cor5}
\end{Theo}

{\it Proof.} For $n=1$, we have done by Corollary \ref{cor3}. Let
$n\ge 2$. Suppose $_AT$ is an $n$-cluster tilting $A$-module such
that $\End_A(T)=\Lambda$, and suppose $_BS$ is an $n$-cluster
tilting $B$-module such that $\End_B(S)= \Gamma$. If $_AM_B$ and
$_BN_A$ are two indecomposable bimodules defining a stable
equivalence of Morita type between $A$ and $B$, then, by the above
discussion, we know that $\Gamma$ is Morita equivalent to
$\End_B(N\otimes_AT)$. Now, we use \cite[Theorem 1.1, or Theorem
1.3]{LX3} which states that if $R$ is an $A$-module with
$\add(_AA)\subseteq \add(R)$, then $\End_A(R)$ and
$\End_B(N\otimes_AR)$ are stably equivalent of Morita type.

Conversely, suppose that two bimodules $_{\Lambda}X_{\Gamma}$ and
$_{\Gamma}Y_{\Lambda}$ define a stable equivalence of Morita type
between $\Lambda$ and $\Gamma$. Note that $\add(T)$ contains both
$\add(_AA)$ and $\add(D(A_A))$. Let $e$ be an idempotent in
$\Lambda$ such that $\add(\Lambda e)$ is just the category of
projective-injective $\Lambda$-modules. Then $e\Lambda e$ is Morita
equivalent to $A^{op}$. As in Theorem \ref{thm}, we have an
idempotent $f$ in $\Gamma$ such that $f\Gamma f$ is Morita
equivalent to $B^{op}$. Thus a stable equivalence of Morita type
between $\Lambda$ and $\Gamma$ implies a stable equivalence of
Morita type between $A^{op}$ and $B^{op}$ by Theorem \ref{thm}, and
therefore a stable equivalence of Morita type between $A$ and $B$.
$\square$

\bigskip
May 20, 2009

\begin{thebibliography}{99}
{\small
\bibitem{ARS}{{\sc M. Auslander, I. Reiten} and {\sc S. O. Smal\o},
Representation thoery of Artin algebras. {\it Cambridge University
Press,} 1995.}

\bibitem{B}{{\sc M. Brou\'e}, Equivalences of blocks of group algebras.
In:{\it Finite dimentional algebras and related topics}. V. Dlab and
L. L. Scott(eds.), Kluwer, (1994), 1-26.}

\bibitem{DM}{{\sc A. Dugas} and {\sc R. Martinez-Villa},
A note on stable equivalences of Morita type . {\it J. Pure Appl.
Algebra} \textbf{208}(2007)421-433.}

\bibitem{hx4}{{\sc W. Hu } and {\sc C. C. Xi}, Derived equivalences of endomorphism and quotient algebras. Preprint,
available at :
http://math.bnu.edu.cn/$^{\sim}$ccxi/Papers/Articles/xihu-4.pdf,
2009.}

\bibitem{hx2}{{\sc W. Hu } and {\sc C. C. Xi}, Almost
$\cal D$-split sequences and derived equivalences. Preprint,
available at :
http://math.bnu.edu.cn/$^{\sim}$ccxi/Papers/Articles/xihu-2.pdf,
2007.}

\bibitem{iyama}{ {\sc O. Iyama}, Cluster tilting for higher Auslander algebras. Preprint, 2008.}

\bibitem{krause}{{\sc H. Krause}, Representation type and stable equivalences of Morita
type for finite dimensional algebras. \emph{Math. Z.}
\textbf{229}(1998)601-606.}

\bibitem{L2}{{\sc M. Linckelmann},
On stable equivalences of Morita type. In:{\it Derived equivalences
for group rings}, LNM \textbf{1685}(1998)221-232.}

\bibitem{LX1}{{\sc Y. M. Liu} and {\sc C. C. Xi},
Constructions of stable equivalences of Morita type for
finite-dimensional algebras. I. {\it Trans. Amer. Math. Soc.} {\bf
358}(2006), no. 6, 2537-2560.}

\bibitem{LX2}{{\sc Y. M. Liu} and {\sc C. C. Xi},
Constructions of stable equivalences of Morita type for
finite-dimensional algebras. II. {\it Math. Z.} {\bf 251}(2005),
no.1, 21-39.}

\bibitem{LX3}{{\sc Y. M. Liu} and {\sc C. C. Xi},
Constructions of stable equivalences of Morita type for
finite-dimensional algebras. III. \emph{J. London Math. Soc.}
\textbf{76}(2007), no. 2, 567-585.}

\bibitem{R1}{{\sc J. Rickard},
The abelian defect group conjecture. {\it Proceedings of the
International Congress of Mathematicians}, Vol. II (Berlin, 1998).
Doc. Math. 1998, Extra Vol. II, 121-128 (electronic).}

\bibitem{R2}{{\sc J. Rickard},
Some recent advances in modular representation theory. {\it Canad.
Math. Soc. Conf. Proc.} {\bf 23}(1998)157-178.}

\bibitem{C3}{{\sc C. C. Xi}, Stable equivalences of adjoint type. {\it Forum Math.} \textbf{20}(2008),
no.1, 81-97.}

\bibitem{C2}{{\sc C. C. Xi}, Representation dimension and
quasi-hereditary algebras. {\it Adv. Math.} {\bf 168}(2002)280-298.}

\bibitem{C5}{{\sc C. C. Xi},
The relative Auslander-Reiten theory of modules. Preprint, available
at: http://math.bnu.edu.cn/$^{\sim}$ccxi/Papers/Articles/rart.pdf,
2005.} }
\end{thebibliography}
\end{document}